# ROBUST NONPARAMETRIC ESTIMATION VIA WAVELET MEDIAN REGRESSION

BY LAWRENCE D. BROWN, T. TONY CAI[1] AND HARRISON H. ZHOU[2]

*University of Pennsylvania, University of Pennsylvania and Yale University*

In this paper we develop a nonparametric regression method that is simultaneously adaptive over a wide range of function classes for the regression function and robust over a large collection of error distributions, including those that are heavy-tailed, and may not even possess variances or means. Our approach is to first use local medians to turn the problem of nonparametric regression with unknown noise distribution into a standard Gaussian regression problem and then apply a wavelet block thresholding procedure to construct an estimator of the regression function. It is shown that the estimator simultaneously attains the optimal rate of convergence over a wide range of the Besov classes, without prior knowledge of the smoothness of the underlying functions or prior knowledge of the error distribution. The estimator also automatically adapts to the local smoothness of the underlying function, and attains the local adaptive minimax rate for estimating functions at a point.

A key technical result in our development is a quantile coupling theorem which gives a tight bound for the quantile coupling between the sample medians and a normal variable. This median coupling inequality may be of independent interest.

**1. Introduction.** A standard nonparametric regression model involves observation of $\{x_i, Y_i\}$ where

(1) $$Y_i = f(x_i) + \xi_i, \qquad i = 1, \ldots, n.$$

Most of the theory that has so far been developed for such a model involves an assumption that the errors $\xi_i$ are independent and identically-distributed

Received June 2007; revised June 2007.
[1]Supported in part by NSF Grants DMS-00-72578 and DMS-03-06576.
[2]Supported in part by NSF Career Award DMS-06-45676.
*AMS 2000 subject classifications.* Primary 62G08; secondary 62G20.
*Key words and phrases.* Adaptivity, asymptotic equivalence, James–Stein estimator, moderate large deviation, nonparametric regression, quantile coupling, robust estimation, wavelets.

(i.i.d.) normal variables. These assumptions are suitable for a wide range of applications of the model. In the Gaussian noise setting many smoothing techniques including wavelet thresholding methods have been developed and shown to be highly adaptive. However, when the noise $\xi_i$ has a heavy-tailed distribution, these techniques are not readily applicable. For example, in Cauchy regression where $\xi_i$ has a Cauchy distribution, typical realizations of $\xi_i$ contain a few extremely large observations of order $n$ since

$$P(\max\{\xi_i\} \geq n) = \left(\frac{1}{\pi}\arctan(n) + \frac{1}{2}\right)^n \to \exp\left(-\frac{1}{\pi}\right).$$

In contrast, the largest observation of the noise $\xi_i$ in Gaussian regression is of order $\sqrt{\log n}$. It is thus clear that the classical denoising methods designed for Gaussian noise would fail if they are applied directly to the sample $\{Y_i\}$ when the noise in fact has a Cauchy distribution. Standard wavelet thresholding procedures would also fail in such a heavy-tailed noise setting. See Section 3.2 for further discussions.

In the usual nonparametric regression case the regression function $f$ is often alternatively described as the conditional expectation $f(x_i) = E(Y_i|x_i)$. However, if the error distributions fail to have a mean, then this conditional expectation will not exist. Even when the conditional expectation exists, estimating the conditional expectation may be a very non-robust goal, and not suitable for particular applications. For error distributions that may be heavy tailed it seems more suitable to estimate the conditional median of $Y_i$. Hence, in the sequel we assume (1) holds with

(2) $\quad\quad\quad\quad\quad\quad \xi_i \text{ i.i.d.} \quad \text{and} \quad \text{median}(\xi_i) = 0.$

There are practical situations for which the normality assumption is not satisfactory. See, for example, Stuck and Kleiner (1974), Stuck (2000) and references therein. It is necessary to develop methods to be used in such cases, and to establish the theoretical properties of these methods. In this paper we develop an estimation method that is simultaneously adaptive over a wide range of function classes for $f$ and robust over a large collection of error distributions for $\xi_i$, including those that are heavy-tailed, and may not even possess variances or means. In brief, our method may be summarized as a blockwise wavelet thresholding implementation built from the medians of suitably binned data. We first divide the interval $[0,1]$ into a number of equal-length subintervals, then take the median of the observations in each subinterval, and finally apply the BlockJS wavelet thresholding procedure developed in Cai (1999) to the local medians together with a bias correction to obtain an estimator of the regression function $f$.

Unlike most wavelet methods, the performance of the algorithm here is not sensitive to the tail behavior of the distribution of $\xi_i$, and hence can be shown to have the necessary robustness property. We show that the estimator enjoys a high degree of adaptivity and robustness. It is shown that the estimator simultaneously attains the exact optimal rate of convergence over a wide range of the Besov classes, without prior knowledge of the smoothness of the underlying function or prior knowledge of the error distribution. The estimator also automatically adapts to the local smoothness of the underlying function, and attains the local adaptive minimax rate for estimating functions at a point.

Donoho and Yu (2000) considered this model for $\alpha$-stable noise, but the risk properties of their proposal are unclear. In the wavelet regression setting, Hall and Patil (1996) studied nonparametric location models and achieved the optimal minimax rate up to a logarithmic term, but under an assumption that $\xi_i$ has a finite fourth moment. As we noted, our results do not need the existence of the mean for the noise or prior knowledge of the error distribution. Most closely related to our work is Averkamp and Houdré (2003, 2005) where the optimal minimax rate of global risk is studied. But their noise is assumed to be known, and their results are not adaptive.

The key technical result in our development is a quantile coupling theorem that is used to connect our problem with a more familiar Gaussian setting. The theorem gives a tight bound for the quantile coupling between the medians of i.i.d. random variables and a normal variable. The result enables us to treat the medians of the observations in the subintervals as if they were normal random variables. The coupling theorem may be of independent interest, since analogous coupling theorems for means have proved to be an important general tool in many contexts. See Section 2 for this result and for further discussion and citations to the literature on quantile coupling.

The paper is organized as follows. In Section 2 we derive a quantile coupling inequality for medians and obtain a moderate large deviation result. This coupling inequality is needed for the proof of the asymptotic properties of our estimation procedure, and may be of independent interest for other statistical applications. Our procedure is defined in Section 3.2 and its asymptotic properties are described in Section 4. Section 5 contains further discussion of our results, and formal proofs are contained in Section 6. The reader interested only in the definition of our wavelet regression procedure and a description of its properties can skip Section 2 and proceed directly to Section 3.

**2. Quantile coupling for median.** We begin with a brief introduction to quantile coupling. Let $X$ be a random variable with distribution $G$ and $Y$ with a continuous distribution $F$. Define

$$\widetilde{X} = G^{-1}(F(Y)), \tag{3}$$

where $G^{-1}(x) = \inf\{u : G(u) \geq x\}$, then $\mathcal{L}(\widetilde{X}) = \mathcal{L}(X)$ [cf. Pollard (2001), page 41]. Note that $\widetilde{X}$ and $Y$ are now defined on the same probability space. This makes it possible to give a pointwise bound between $\widetilde{X}$ and $Y$. The first tight bound of quantile coupling between the sum of i.i.d. random variables with a normal random variable was given in Komlós, Major and Tusnády (1975). A bound for the coupling of a Binomial random variable with a normal random variable is given as follows. For $X \sim \text{Binomial}(n, 1/2)$ and $Y \sim N(n/2, n/4)$, let $\widetilde{X}(Y)$ be defined as in equation (3). Then for some constant $C > 0$ and $\varepsilon > 0$, when $|\widetilde{X}| \leq \varepsilon n$,

$$|\widetilde{X} - Y| \leq C + C \frac{|\widetilde{X}|^2}{n}. \tag{4}$$

This result plays a key role in the KMT/Hungarian construction to couple the empirical distribution with a Brownian bridge. A detailed proof of the result can be found in Mason (2001) and Bretagnolle and Massart (1989). A general theory for improving the classical quantile coupling bound was given in Zhou (2005).

Standard coupling inequalities are mostly focused on the coupling of the mean of i.i.d. random variables with a normal variable. In this section we study the coupling of a median statistic with a normal variable. We derive a moderate deviation result for the median statistic and obtain a quantile coupling inequality similar to the classical KMT bound for the mean. This coupling result plays a crucial role in this paper. It is the main tool for reducing the problem of robust estimation with unknown noise to a well studied problem of Gaussian regression with unknown variance. The result here may be of independent interest because of the fundamental role played by the median in statistics.

Let $X_1, \ldots, X_n$ be i.i.d. random variables with density function $h$. Denote the sample median by $X_{\text{med}}$. We will construct a new random variable $\widetilde{X}_{\text{med}}$ by using quantile coupling in (3) such that $\mathcal{L}(\widetilde{X}_{\text{med}}) = \mathcal{L}(X_{\text{med}})$ and show that $\widetilde{X}_{\text{med}}$ can be well approximated by a normal random variable as equation (4). We need the following assumptions on the density function $h(x)$ to derive the quantile coupling inequality.

ASSUMPTION (A1). $\int_{-\infty}^{0} h(x) = \frac{1}{2}$, $h(0) > 0$, and $h(x)$ is Lipschitz at $x = 0$.

Here the Lipschitz condition at 0 means that there is a constant $C > 0$ such that $|h(x) - h(0)| \leq C|x|$ in an open neighborhood of 0. This condition implies that $h$ is continuous at 0. We assume $h(0) > 0$ so that the median of the distribution is unique and the distribution of the sample median is asymptotically normal [cf. Casella and Berger (2002), page 483]. The Lipschitz condition is assumed so that a moderate large deviation result for the distribution of sample median can be obtained to derive a quantile coupling inequality as in equation (4).

THEOREM 1. *Let $Z$ be a standard normal random variable and let $X_1, \ldots, X_n$ be i.i.d. with density function $h$ where $n = 2k+1$ for some integer $k \geq 1$. Let Assumption* (A1) *hold. Then for every $n$ there is a mapping $\widetilde{X}_{\mathrm{med}}(Z) : \mathbb{R} \mapsto \mathbb{R}$ such that $\mathcal{L}(\widetilde{X}_{\mathrm{med}}(Z)) = \mathcal{L}(X_{\mathrm{med}})$ and*

(5) $\quad |\sqrt{4n} h(0) \widetilde{X}_{\mathrm{med}} - Z| \leq \dfrac{C}{\sqrt{n}} + \dfrac{C}{\sqrt{n}} |\sqrt{4n} h(0) \widetilde{X}_{\mathrm{med}}|^2 \qquad$ *when $|\widetilde{X}_{\mathrm{med}}| \leq \varepsilon$*

*where $C, \varepsilon > 0$ depend on $h$ but not on $n$.*

The quantile coupling bound here is similar to the classical KMT bound (4) for the sample mean. This result has close connection to strong approximation of quantile process in Csörgő and Révész (1978). The condition of Theorem 1 here is weaker. Only a Lipschitz condition at $x = 0$ is assumed here to establish the non-uniform bound given in (5). As shown in Zhou (2005), the classical quantile coupling bound for the mean can be improved when the distribution of $X_i$ is symmetric. Similarly, if we assume $h'(0) = 0$, the bound in Theorem 1 can be improved from the rate $1/\sqrt{n}$ to the rate $1/n$. See section 4 for more details. The bound in Theorem 1 can also be expressed in terms of $Z$, as follows.

COROLLARY 1. *Under the assumption of Theorem 1, the mapping $\widetilde{X}_{\mathrm{med}}(Z)$ in Theorem 1 satisfies*

(6) $\qquad |\sqrt{4n} h(0) \widetilde{X}_{\mathrm{med}} - Z| \leq \dfrac{C}{\sqrt{n}} (1 + |Z|^2) \qquad$ *when $|Z| \leq \varepsilon \sqrt{n}$*

*where $C, \varepsilon > 0$ do not depend on $n$.*

REMARK 1. When $n = 2k$ is even, the sample median $X_{\mathrm{med}}$ is usually taken to be $(X_{(k)} + X_{(k+1)})/2$. Similar quantile coupling inequalities as (5) and (6) can be obtained. For each $i$, let $X_{-i,\mathrm{med}}$ be the median of the original sample with $X_i$ removed. Then $X_{\mathrm{med}} = \frac{1}{n} \sum_{i=1}^{n} X_{-i,\mathrm{med}}$. Let $G_{n-1}$ be the distribution of the median of $n-1$ i.i.d. observations with density $h$ and define $(Z_i)_{1 \leq i \leq n} \sim \mathcal{L}(\Phi^{-1} \circ G_{n-1}(X_{-i,\mathrm{med}}), 1 \leq i \leq n)$. Let $\widetilde{X}_{-i,\mathrm{med}} = G_{n-1}^{-1} \Phi(Z_i)$. Then $\mathcal{L}(\widetilde{X}_{-i,\mathrm{med}}, 1 \leq i \leq n) = \mathcal{L}(X_{-i,\mathrm{med}}, 1 \leq i \leq n)$. Now a direct application of Theorem 1 gives

$$|\widetilde{X}_{\mathrm{med}} - Z| \leq \dfrac{C}{\sqrt{n}} (1 + |\sqrt{4n} h(0) (|\widetilde{X}_{(k)}| + |\widetilde{X}_{(k+1)}|)|^2)$$

when $|\widetilde{X}_{(k)}| + |\widetilde{X}_{(k+1)}| \leq \varepsilon$, and $Z = \frac{1}{n} \sum_{i=1}^{n} Z_i$. So in Sections 3 and 5 we assume the number of observations in each bin is odd without loss of generality.

The coupling result given in Theorem 1 in fact holds uniformly over a rich collection of distributions. For $0 < \epsilon_1 < 1$ and $\epsilon_2 > 0$ define

$$\mathcal{H}_{\epsilon_1,\epsilon_2} = \bigg\{ h \colon \int_{-\infty}^{0} h(x) = \frac{1}{2}, \epsilon_1 \leq h(0) \leq \frac{1}{\epsilon_1},$$

(7)

$$|h(x) - h(0)| \leq \frac{|x|}{\epsilon_1} \text{ for all } |x| < \epsilon_2 \bigg\}.$$

It can be shown that Theorem 1 holds uniformly for the whole family of $h \in \mathcal{H}_{\epsilon_1,\epsilon_2}$.

THEOREM 2. *Let $X_1, \ldots, X_n$ be i.i.d. with density $h \in \mathcal{H}_{\epsilon_1,\epsilon_2}$. For every $n = 2k+1$ with integer $k \geq 1$, there is a mapping $\widetilde{X}_{\mathrm{med}}(Z) \colon \mathbb{R} \mapsto \mathbb{R}$ such that $\mathcal{L}(\widetilde{X}_{\mathrm{med}}(Z)) = \mathcal{L}(X_{\mathrm{med}})$ and for two constants $C_{\epsilon_1,\epsilon_2}$, $\varepsilon_{\epsilon_1,\epsilon_2} > 0$ depending only on $\epsilon_1$ and $\epsilon_2$*

$$|\sqrt{4n} h(0) \widetilde{X}_{\mathrm{med}} - Z| \leq \frac{C_{\epsilon_1,\epsilon_2}}{\sqrt{n}} + \frac{C_{\epsilon_1,\epsilon_2}}{\sqrt{n}} |\sqrt{4n} h(0) \widetilde{X}_{\mathrm{med}}|^2$$

*uniformly over all $h \in \mathcal{H}_{\epsilon_1,\epsilon_2}$.*

REMARK 2. The quantile coupling inequalities in Corollary 1 and Remark 1 also hold uniformly over $\mathcal{H}_{\epsilon_1,\epsilon_2}$ by replacing $C$ and $\varepsilon$ there with two constants depending $\epsilon_1$ and $\epsilon_2$.

**3. Methodology for robust wavelet regression.** We now define our robust nonparametric regression estimator. Then we apply the median quantile coupling results developed in the previous section to establish its asymptotic properties.

As we have mentioned, the first key step in our approach is to bin the data according to the values of the independent variable. The sample median is then computed within each bin. This leads to a new data situation in which the bin centers are treated as the independent variables in a nonparametric regression, with the bin medians being the dependent variables. This new situation can then be satisfactorily viewed as if it were a Gaussian regression problem. It is important that the number of bins be chosen in a suitable range. For the applications in our paper it turns out to be appropriate to choose the number of bins to be $T \asymp n^{3/4}$, where $n$ is the original sample size. It appears that such a choice of $T$ would also be suitable for use with many other Gaussian nonparametric regression methods.

Proceeding in this way one should expect as a heuristic principle that the resulting nonparametric procedure will inherit the asymptotic optimality properties of the Gaussian nonparametric regression technique that is employed. Of course, this heuristic principle needs to be established in particular cases. The difficulty of doing so will depend on the nature of the Gaussian technique and the generality of the asymptotic assumptions.

In the present treatment we choose to employ a Gaussian wavelet method involving a block James–Stein wavelet estimator. Implementation of the procedure is straightforward since the number of bins can be chosen as a power of 2, as is especially convenient for wavelet implementation. This estimator enjoys excellent asymptotic adaptivity properties in the Gaussian setting. We show that the current binned-median version has analogous properties over nearly the same range of Besov balls as does the original Gaussian procedure. The precise statement of asymptotic properties is contained in Theorems 3 and 4. The full strength of the asymptotic properties of our wavelet procedure in a Gaussian setting depends on detailed moderate-deviation properties of the Gaussian distribution. For this reason our proof of asymptotic properties of the binned median version requires careful treatment of moderate-deviation properties of the binned medians, as in the coupling results established in Section 2.

We shall focus on the case where the design points $\{x_i\}$, are equally spaced on the interval $[0, 1]$. The more general case will be discussed at the end of Section 4. The procedure, which will be described in detail in the next section, can be briefly summarized as follows. Let the sample $\{Y_i, i = 1, \ldots, n\}$ be given as in (1) where $x_i = \frac{i}{n}$ and the noise variables $\xi_i$ are i.i.d. with an unknown density $h$. Let $J = \lfloor \log_2 n^{3/4} \rfloor$. Set $T = 2^J$ and $m = n/T$. We divide the interval $[0, 1]$ into $T$ equal-length subintervals. Note that $T \asymp n^{3/4}$. For $1 \leq j \leq T$, let $I_j = \{Y_i : x_i \in (\frac{j-1}{T}, \frac{j}{T}]\}$ be the $j$th bin and let $X_j$ be the median of the observations in $I_j$. We treat $X_j$ as if it were a normal random variable with mean $f(\frac{j}{T}) + b_m$ and variance $1/(4mh^2(0))$ (see Theorem 1), where

(9) $$b_m = E\{\text{median}(\xi_1, \ldots, \xi_m)\}.$$

Then apply a nonparametric Gaussian regression procedure. In this paper, we apply the BlockJS wavelet thresholding procedure developed in Cai (1999) to construct an estimator of $f$. The final estimator $\hat{f}$ is given in equations (16) and (18).

We begin in Section 3.1 with a brief introduction to wavelet block thresholding in the Gaussian regression setting and then give a detailed description of our wavelet procedure for robust estimation in Section 3.2.

3.1. *Wavelet block thresholding for Gaussian regression.* Let $\{\phi, \psi\}$ be a pair of father and mother wavelets. The functions $\phi$ and $\psi$ are assumed to be compactly supported and $\int \phi = 1$. Dilation and translation of $\phi$ and $\psi$ generates an orthonormal wavelet basis. For simplicity in exposition, in the present paper we work with periodized wavelet bases on $[0, 1]$. Let

$$\phi_{j,k}^p(t) = \sum_{l=-\infty}^{\infty} \phi_{j,k}(t-l), \qquad \psi_{j,k}^p(t) = \sum_{l=-\infty}^{\infty} \psi_{j,k}(t-l) \qquad \text{for } t \in [0,1]$$

where $\phi_{j,k}(t) = 2^{j/2}\phi(2^j t - k)$ and $\psi_{j,k}(t) = 2^{j/2}\psi(2^j t - k)$. The collection $\{\phi_{j_0,k}^p, k = 1, \ldots, 2^{j_0}; \psi_{j,k}^p, j \geq j_0 \geq 0, k = 1, \ldots, 2^j\}$ is then an orthonormal basis of $L^2[0,1]$, provided the primary resolution level $j_0$ is large enough to ensure that the support of the scaling functions and wavelets at level $j_0$ is not the whole of $[0,1]$. The superscript "$p$" will be suppressed from the notation for convenience. An orthonormal wavelet basis has an associated orthogonal Discrete Wavelet Transform (DWT) which transforms sampled data into the wavelet coefficients. See Daubechies (1992) and Strang (1992) for further details about the wavelets and discrete wavelet transform. A square-integrable function $f$ on $[0, 1]$ can be expanded into a wavelet series:

$$(10) \qquad f(t) = \sum_{k=1}^{2^{j_0}} \tilde{\theta}_{j_0,k}\phi_{j_0,k}(t) + \sum_{j=j_0}^{\infty} \sum_{k=1}^{2^j} \theta_{j,k}\psi_{j,k}(t)$$

where $\tilde{\theta}_{j,k} = \langle f, \phi_{j,k} \rangle, \theta_{j,k} = \langle f, \psi_{j,k} \rangle$ are the wavelet coefficients of $f$.

The BlockJS procedure was proposed in Cai (1999) for Gaussian nonparametric regression and was shown to achieve simultaneously three objectives: adaptivity, spatial adaptivity, and computational efficiency. The procedure can be most easily explained in the sequence space setting. Suppose we observe the wavelet sequence data:

$$(11) \qquad y_{j,k} = \theta_{j,k} + \sigma z_{j,k}, \qquad j \geq j_0, k = 1, 2, \ldots, 2^j$$

where $z_{j,k}$ are i.i.d. $N(0,1)$ and the noise level $\sigma$ is known. The mean vector $\theta$ is the object of interest. The BlockJS procedure is as follows. Let $J = [\log_2 n]$. Divide each resolution level $j_0 \leq j < J$ into nonoverlapping blocks of length $L = [\log n]$ (or $L = 2^{\lfloor \log_2(\log n) \rfloor} \approx \log n$). Let $B_j^i$ denote the set of indices of the coefficients in the $i$-th block at level $j$, that is, $B_j^i = \{(j,k) : (i-1)L + 1 \leq k \leq iL\}$. Let $S_{j,i}^2 \equiv \sum_{(j,k) \in B_j^i} y_{j,k}^2$ denote the sum of squared empirical wavelet coefficients in block $B_j^i$. A James–Stein type shrinkage rule is then applied to each block $B_j^i$. For $(j, k) \in B_j^i$,

$$(12) \qquad \hat{\theta}_{j,k} = \begin{cases} \left(1 - \dfrac{\lambda_* L \sigma^2}{S_{j,i}^2}\right)_+ y_{j,k}, & \text{for } (j,k) \in B_j^i, \ j_0 \leq j < J, \\ 0, & \text{for } j \geq J, \end{cases}$$

where $\lambda_* = 4.50524$ is a constant satisfying $\lambda_* - \log \lambda_* = 3$. The threshold $\lambda_* = 4.50524$ is selected according to a block thresholding oracle inequality and a minimax criterion. See Cai (1999) for further details.

3.2. *Wavelet procedure for robust regression.* Now we are ready to give a detailed description of our procedure for robust estimation. Hereafter we shall set $g(t) = f(t) + b_m$ where $b_m$ is given as in (9).

Apply the discrete wavelet transform to the binned medians $X = (X_1, \ldots, X_T)$, and let $U = T^{-1/2}WX$ be the empirical wavelet coefficients, where $W$ is the discrete wavelet transformation matrix. Write

$$(13) \quad U = (\widetilde{y}_{j_0,1}, \ldots, \widetilde{y}_{j_0,2^{j_0}}, y_{j_0,1}, \ldots, y_{j_0,2^{j_0}}, \ldots, y_{J-1,1}, \ldots, y_{J-1,2^{J-1}})'.$$

Here $\widetilde{y}_{j_0,k}$ are the gross structure terms at the lowest resolution level, and $y_{j,k}$ $(j = j_0, \ldots, J-1, k = 1, \ldots, 2^j)$ are empirical wavelet coefficients at level $j$ which represent fine structure at scale $2^j$. The empirical wavelet coefficients can be written as

$$(14) \quad y_{j,k} = \theta_{j,k} + \epsilon_{j,k} + \frac{1}{2h(0)\sqrt{n}} z_{j,k} + \xi_{j,k},$$

where $\theta_{j,k}$ are the true wavelet coefficients of $g = f + b_m$, $\epsilon_{j,k}$ are "small" deterministic approximation errors, $z_{j,k}$ are i.i.d. $N(0,1)$, and $\xi_{j,k}$ are some "small" stochastic errors. The theoretical calculations given in Section 6 will show that both the approximation errors $\epsilon_{j,k}$ and the stochastic errors $\xi_{j,k}$ are negligible in certain sense. If these negligible errors are ignored then we have

$$(15) \quad y_{j,k} \approx \theta_{j,k} + \frac{1}{2h(0)\sqrt{n}} z_{j,k},$$

which is the same as the idealized sequence model (11) with noise level $\sigma = 1/(2h(0)\sqrt{n})$.

The BlockJS procedure is then applied to the empirical coefficients $y_{j,k}$ as if they are distributed as in (15). More specifically, at each resolution level $j$, the empirical wavelet coefficients $y_{j,k}$ are grouped into nonoverlapping blocks of length $L$. As in the sequence estimation setting let $B_j^i = \{(j,k) : (i-1)L + 1 \leq k \leq iL\}$ and let $S_{j,i}^2 \equiv \sum_{(j,k) \in B_j^i} y_{j,k}^2$. Let $\widehat{h}^2(0)$ be an estimator of $h^2(0)$ [see equation (38) for an estimator]. A modified James–Stein shrinkage rule is then applied to each block $B_j^i$, that is,

$$(16) \qquad \hat{\theta}_{j,k} = \left(1 - \frac{\lambda_* L}{4\widehat{h}^2(0)nS_{j,i}^2}\right)_+ y_{j,k} \qquad \text{for } (j,k) \in B_j^i,$$

where $\lambda_* = 4.50524$ is the solution to the equation $\lambda_* - \log \lambda_* = 3$ and $4\widehat{h}^2(0)n$ in the shrinkage factor of (16) is due to the fact that the noise level in (15) is $\sigma = 1/(2h(0)\sqrt{n})$. For the gross structure terms at the lowest resolution level $j_0$, we set $\hat{\widetilde{\theta}}_{j_0,k} = \widetilde{y}_{j_0,k}$. The estimate of $g$ at the equally spaced sample points $\{\frac{i}{T} : i = 1, \ldots, T\}$ is then obtained by applying the inverse discrete wavelet transform (IDWT) to the denoised wavelet coefficients. That is, $\{g(\frac{i}{T}): i=1,\ldots,T\}$ is estimated by $\widehat{g} = \{\widehat{g(\frac{i}{T})}: i=1,\ldots,T\}$ with $\widehat{g} = T^{1/2}W^{-1} \cdot \hat{\theta}$. The estimate of the whole function $g = f + b_m$ is given by

$$\widehat{g}(t) = \sum_{k=1}^{2^{j_0}} \hat{\tilde{\theta}}_{j_0,k}\phi_{j_0,k}(t) + \sum_{j=j_0}^{J-1}\sum_{k=1}^{2^j} \hat{\theta}_{j,k}\psi_{j,k}(t).$$

To get an estimator of $f$ we need to also estimate $b_m$. This is done as follows. Divide each bin $I_j$ into two sub-bins with the first bin of the size $\lfloor \frac{m}{2} \rfloor$. Let $X_j^*$ be the median of observations in the first sub-bin. We set

$$(17) \qquad \hat{b}_m = \frac{1}{T}\sum_j (X_j^* - X_j)$$

and define

$$(18) \quad \widehat{f}_n(t) = \widehat{g}_n(t) - \hat{b}_m = \sum_{k=1}^{2^{j_0}} \hat{\tilde{\theta}}_{j_0,k}\phi_{j_0,k}(t) + \sum_{j=j_0}^{J-1}\sum_{k=1}^{2^j} \hat{\theta}_{j,k}\psi_{j,k}(t) - \hat{b}_m.$$

REMARK 3. The quantity $b_m$ is the systematic bias due to the expectation of the median of the noise $\xi_i$ in each bin. Lemma 5 in Section 6 shows that $b_m = -\frac{h'(0)}{8h^3(0)}m^{-1} + O(m^{-2})$. Hence this systematic bias can possibly be dominant if it is ignored. The estimate $\hat{b}_m$ serves as "bias correction." Lemma 5 shows that the estimation error of $\hat{b}_m$ is negligible relative to the minimax risk of $\hat{f}_n$ when $m = O(n^{1/4})$.

**4. Adaptivity and robustness of the procedure.** We study the theoretical properties of our procedure over the Besov spaces that are by now standard for the analysis of wavelet regression methods. Besov spaces are a very rich class of function spaces and contain as special cases many traditional smoothness spaces such as Hölder and Sobolev spaces. Roughly speaking, the Besov space $B_{p,q}^{\alpha}$ contains functions having $\alpha$ bounded derivatives in $L^p$ norm, the third parameter $q$ gives a finer gradation of smoothness. Full details of Besov spaces are given, for example, in Triebel (1992) and DeVore and Popov (1988). For a given $r$-regular mother wavelet $\psi$ with $r > \alpha$ and a fixed primary resolution level $j_0$, the Besov sequence norm $\|\cdot\|_{b_{p,q}^{\alpha}}$ of the wavelet coefficients of a function $f$ is then defined by

$$(19) \qquad \|f\|_{b_{p,q}^{\alpha}} = \|\underline{\xi}_{j_0}\|_p + \left(\sum_{j=j_0}^{\infty}(2^{js}\|\underline{\theta}_j\|_p)^q\right)^{1/q}$$

where $\underline{\xi}_{j_0}$ is the vector of the father wavelet coefficients at the primary resolution level $j_0$, $\underline{\theta}_j$ is the vector of the wavelet coefficients at level $j$, and $s = \alpha + \frac{1}{2} - \frac{1}{p} > 0$. Note that the Besov function norm of index $(\alpha, p, q)$ of a function $f$ is equivalent to the sequence norm (19) of the wavelet coefficients of the function. See Meyer (1992). We define

$$(20) \quad B_{p,q}^\alpha(M) = \{f; \|f\|_{b_{p,q}^\alpha} \leq M\}.$$

In the case of Gaussian noise Donoho and Johnstone (1998) show that the minimax risk of estimating $f$ over the Besov body $B_{p,q}^\alpha(M)$,

$$(21) \quad R^*(B_{p,q}^\alpha(M)) = \inf_{\hat{f}} \sup_{f \in B_{p,q}^\alpha(M)} E\|\hat{f} - f\|_2^2,$$

converges to 0 at the rate of $n^{-2\alpha/(1+2\alpha)}$ as $n \to \infty$.

In addition to Assumption (A1) in Section 2, we need the following weak condition on the density $h$ of $\xi_i$.

ASSUMPTION (A2). $\int |x|^{\epsilon_3} h(x) \, dx < \infty$ for some $\epsilon_3 > 0$.

This assumption guarantees that the moments of the median of the binned data are well approximated by those of the normal random variable. Note that Assumption (A2) is satisfied by Cauchy distribution for any $0 < \epsilon_3 < 1$. For $0 < \epsilon_1 < 1$, $\epsilon_i > 0$, $i = 2, 3, 4$, define $\mathcal{H} = \mathcal{H}(\epsilon_1, \epsilon_2, \epsilon_3, \epsilon_4)$ by

$$(22) \quad \mathcal{H} = \left\{ h : h \in \mathcal{H}_{\epsilon_1, \epsilon_2}, |h^{(3)}(x)| \leq \epsilon_4 \text{ for } |x| \leq \epsilon_3 \text{ and } \int |x|^{\epsilon_3} h(x) \, dx < \epsilon_4 \right\}.$$

The following theorem shows that our estimator achieves optimal global adaptation for a wide range of Besov balls $B_{p,q}^\alpha(M)$ defined in (20) and uniformly over the family of noise distributions given in (22).

THEOREM 3. *Suppose the wavelet $\psi$ is $r$-regular. Then the estimator $\hat{f}_n$ defined in (18) satisfies, for $p \geq 2$, $\alpha \leq r$ and $\frac{2\alpha^2 - \alpha/3}{1 + 2\alpha} > \frac{1}{p}$,*

$$\sup_{h \in \mathcal{H}} \sup_{f \in B_{p,q}^\alpha(M)} E\|\hat{f}_n - f\|_2^2 \leq C n^{-2\alpha/(1+2\alpha)},$$

*and for $1 \leq p < 2$, $\alpha \leq r$ and $\frac{2\alpha^2 - \alpha/3}{1 + 2\alpha} > \frac{1}{p}$,*

$$\sup_{h \in \mathcal{H}} \sup_{f \in B_{p,q}^\alpha(M)} E\|\hat{f}_n - f\|_2^2 \leq C n^{-2\alpha/(1+2\alpha)} (\log n)^{(2-p)/(p(1+2\alpha))}.$$

Theorem 3 shows that the estimator simultaneously attains the optimal rate of convergence over a wide range of the Besov classes for $f$ and a large collection of the unknown error distributions for $\xi_i$. In this sense, the estimator enjoys a high degree of adaptivity and robustness.

For functions of spatial inhomogeneity, the local smoothness of the functions varies significantly from point to point and global risk given in Theorem 3 cannot wholly reflect the performance of estimators at a point. The local risk measure

$$R(\widehat{f}(t_0), f(t_0)) = E(\widehat{f}(t_0) - f(t_0))^2 \tag{23}$$

is used for spatial adaptivity.

The local smoothness of a function can be measured by its local Hölder smoothness index. For a fixed point $t_0 \in (0,1)$ and $0 < \alpha \leq 1$, define the local Hölder class $\Lambda^\alpha(M, t_0, \delta)$ as follows:

$$\Lambda^\alpha(M, t_0, \delta) = \{f : |f(t) - f(t_0)| \leq M|t - t_0|^\alpha, \text{ for } t \in (t_0 - \delta, t_0 + \delta)\}.$$

If $\alpha > 1$, then

$$\Lambda^\alpha(M, t_0, \delta) = \{f : |f^{(\lfloor \alpha \rfloor)}(t) - f^{(\lfloor \alpha \rfloor)}(t_0)| \leq M|t - t_0|^{\alpha'} \text{ for } t \in (t_0 - \delta, t_0 + \delta)\}$$

where $\lfloor \alpha \rfloor$ is the largest integer less than $\alpha$ and $\alpha' = \alpha - \lfloor \alpha \rfloor$.

In Gaussian nonparametric regression setting, it is a well-known fact that for estimation at a point, one must pay a price for adaptation. The optimal rate of convergence for estimating $f(t_0)$ over function class $\Lambda^\alpha(M, t_0, \delta)$ with $\alpha$ completely known is $n^{-2\alpha/(1+2\alpha)}$. Lepski (1990) and Brown and Low (1996a, 1996b) showed that one has to pay a price for adaptation of at least a logarithmic factor. It is shown that the local adaptive minimax rate over the Hölder class $\Lambda^\alpha(M, t_0, \delta)$ is $(\log n/n)^{2\alpha/(1+2\alpha)}$.

The following theorem shows that our estimator achieves optimal local adaptation with the minimal cost uniformly over the family of noise distributions defined in (22).

THEOREM 4. *Suppose the wavelet $\psi$ is $r$-regular with $r \geq \alpha > 1/6$. Let $t_0 \in (0,1)$ be fixed. Then the estimator $\hat{f}_n$ defined in (18) satisfies*

$$\sup_{h \in \mathcal{H}} \sup_{f \in \Lambda^\alpha(M, t_0, \delta)} E(\widehat{f}_n(t_0) - f(t_0))^2 \leq C \cdot \left(\frac{\log n}{n}\right)^{2\alpha/(1+2\alpha)}. \tag{24}$$

Theorem 4 shows that the estimator automatically attains the local adaptive minimax rate for estimating functions at a point, without prior knowledge of the smoothness of the underlying functions or prior knowledge of the error distribution.

REMARK 4. After binning and taking the medians, in principle any standard wavelet thresholding estimators could then be used. For example, the VisuShrink procedure of Donoho and Johnstone (1994) with threshold $\lambda = \sigma\sqrt{2\log n}$ can be applied. In this case the resulting estimator satisfies

$$\sup_{h \in \mathcal{H}} \sup_{f \in B^\alpha_{p,q}(M)} E\|\widehat{f}_n - f\|_2^2 \leq C \left(\frac{\log n}{n}\right)^{2\alpha/(1+2\alpha)}$$

for $1 \leq p \leq \infty$, $\alpha \leq r$ and $\frac{2\alpha^2 - \alpha/3}{1+2\alpha} > \frac{1}{p}$ and

(25) $$\sup_{h \in \mathcal{H}} \sup_{f \in \Lambda^\alpha(M, t_0, \delta)} E(\widehat{f}_n(t_0) - f(t_0))^2 \leq C \cdot \left(\frac{\log n}{n}\right)^{2\alpha/(1+2\alpha)}$$

for $r \geq \alpha > 1/6$.

We have so far focused on the equally spaced design case. When the design is not equally spaced, one can either group the sample using equal-length subintervals as in Section 3.2 or bin the sample so that each bin contains the same number of observations, and then take the median of each bin. The first method produces equally spaced medians that are heteroskedastic with the variances depending on the number of observations in the bins. In this case a wavelet procedure for heteroskedastic Gaussian noise can then be applied to the medians to obtain an estimator of $f$. The second method produces unequally spaced medians that are homoskedastic since the number of observations in the bins are the same. A wavelet procedure for unequally spaced observations with homoskedastic Gaussian noise can then be used to get an estimator of $f$. For wavelet procedures for heteroskedastic Gaussian noise or unequally spaced samples, see, for example, Cai and Brown (1998), Kovac and Silverman (2000) and Antoniadis and Fan (2001).

**5. Further discussion.** Theorem 1 gives a general quantile coupling inequality between the median of i.i.d. random variables $X_1, \ldots, X_n$ and a normal random variable. The collection of the distributions of the i.i.d. random variables includes the Cauchy and Gaussian distributions as special cases. Note that for both Cauchy and Gaussian distributions, $h'(0) = 0$, which suggests we may have a tighter quantile coupling bound as in Zhou (2005). Let us further assume that $h'(0) = 0$, and $h''(0)$ exists. We can derive a sharper moderate large deviation result for the median and then obtain a tighter quantile coupling inequality which improves the classical quantile coupling bounds with a rate $1/\sqrt{n}$ under certain smoothness conditions for the distribution function. For every $n$, we can show that there is a mapping $\widetilde{X}_{\mathrm{med}}(Z): \mathbb{R} \mapsto \mathbb{R}$ such that the random variable $\widetilde{X}_{\mathrm{med}}(Z)$ has the same distribution as the median $X_{\mathrm{med}}$ of $X_1, \ldots, X_n$ and

$$|\sqrt{4n}h(0)\widetilde{X}_{\mathrm{med}} - Z| \leq C\frac{1}{n}(1 + |Z|^3) \quad \text{when } |Z| \leq \varepsilon\sqrt{n}$$

where $C, \varepsilon > 0$ do not depend on $n$. We can even establish an asymptotic equivalence result in Le Cam's sense. Assume that

$$f \in \mathcal{F} = \{f : |f(y) - f(x)| \leq M|x - y|^d\}$$

with $d > 3/4$. In the current setting, we modify the procedure with $T = n^{2/3}/\log n$. Then $m = n/T = n^{1/3} \log n$. Recall that $X_j$ is the median of the observations on each bin $I_j$ with $1 \leq j \leq T$. Let $\eta_j$ be the median of corresponding noise, then

$$\min_{(j-1)m+1 \leq i \leq jm} f\left(\frac{i}{n}\right) \leq X_j - \eta_j \leq \max_{(j-1)m+1 \leq i \leq jm} f\left(\frac{i}{n}\right).$$

We need to give an asymptotic justification that it is fine treating $X_j$ as if it were a normal random variable with mean $f(j/T)$ and variance $\frac{1}{4mh^2(0)}$. We can show that observing $\{X_j\}$ is asymptotically equivalent to observing

$$X_j^\dagger = f\left(\frac{j}{T}\right) + Z_j, \qquad Z_j \overset{i.i.d.}{\sim} N\left(0, \frac{1}{4mh^2(0)}\right), \qquad 1 \leq j \leq T$$

in Le Cam's sense by showing that the total variation distance between the distributions of $X_j$'s and $X_j^\dagger$'s tends to 0, that is,

$$|\mathcal{L}(\{X_j\}) - \mathcal{L}(\{X_j^\dagger\})|_{\text{TV}} \to 0.$$

The result shows that asymptotically there is no difference between observing $X_j$'s and observing $X_j^\dagger$'s. That means all optimal statistical procedures for the Gaussian model can be carried over to nonparametric robust estimation for bounded losses. For instance, the asymptotic equivalence here implies that adaptive procedures including SureShrink of Donoho and Johnstone (1995), the empirical Bayes estimation of Zhang (2005) and SureBlock of Cai and Zhou (2006) can be carried over from the Gaussian regression to the Cauchy regression or more general regression. The details of our results will be reported elsewhere. Readers may find recent developments in the asymptotic equivalence theory in Brown and Low (1996a, 1996b), Nussbaum (1996), Grama and Nussbaum (1998) and Golubev, Nussbaum and Zhou (2005).

**6. Proofs.** We shall prove the main results in the order of Theorem 3, Theorem 4 and then Theorems 1 and 2. In this section $C$ denotes a positive constant not depending on $n$ that may vary from place to place and we set $d \equiv \min(\alpha - \frac{1}{p}, 1)$. For simplicity we shall assume that $n$ is divisible by $T$ in the proof. We first collect necessary tools that are needed for the proofs of Theorems 3 and 4.

6.1. *Preparatory results.* In our procedure, there are two steps: (1) binning the data and taking the median in each bin; (2) applying wavelet transform to the medians and using BlockJS to construct an estimator of $f$. In this section, we give two results associated with these two steps. Recall that we denote by $X_j$ the median of each bin $I_j$ in step 1 and treat $X_j$ as if it were a normal random variable with mean $f(j/T) - b_m$ and variance $1/(4mh^2(0))$.

The coupling inequality and the fact that a Besov ball $B_{p,q}^{\alpha}(M)$ can be embedded into a Hölder ball with smoothness $d = \min(\alpha - \frac{1}{p}, 1) > 0$ [cf. Meyer (1992)] enable us to precisely control the difference between $X_j$ and that normal variable. Proposition 1 gives the bounds for both the deterministic and stochastic errors. In Proposition 2 we obtain two risk bounds for the BlockJS procedure used in step 2. These two bounds are used to study global and local adaptation in the following sections.

PROPOSITION 1. *Let $X_j$ be given as in our procedure and let $f \in B_{p,q}^{\alpha}(M)$. Then $X_j$ can be written as*

$$\sqrt{m} X_j = \sqrt{m} f\left(\frac{j}{T}\right) + \sqrt{m} b_m + \frac{1}{2} Z_j + \epsilon_j + \zeta_j \tag{26}$$

*where:*

(i) $Z_j \overset{i.i.d.}{\sim} N(0, \frac{1}{h^2(0)})$;

(ii) $\epsilon_j$ *are constants satisfying* $|\epsilon_j| \leq C\sqrt{m} T^{-d}$ *and so* $\frac{1}{n} \sum_{i=1}^{T} \epsilon_j^2 \leq C T^{-2d}$;

(iii) $\zeta_j$ *are independent and "stochastically small" random variables satisfying with $E\zeta_j = 0$, for any $l > 0$*

$$E|\zeta_j|^l \leq C_l m^{-l/2} + C_l m^{l/2} T^{-dl} \tag{27}$$

*and for any $a > 0$*

$$P(|\zeta_j| > a) \leq C_l (a^2 m)^{-l/2} + C_l (a^2 T^{2d}/m)^{-l/2} \tag{28}$$

*where $C_l > 0$ is a constant depending on $l$ only.*

PROOF. Let $\eta_j = \text{median}(\{\xi_i : (j-1)m + 1 \leq i \leq jm\})$. We define $Z_j = \frac{1}{h(0)} \Phi^{-1}(G(\eta_j))$ where $G$ is the distribution of $\eta_j$. It follows from Theorem 1 that $\sqrt{4m}\eta_j$ is well approximated by $Z_j$ whose distribution is $N(0, \frac{1}{h^2(0)})$. Set

$$\epsilon_j = \sqrt{m} E X_j - \sqrt{m} f\left(\frac{j}{T}\right) - \sqrt{m} b_m$$
$$= E\left\{\sqrt{m} X_j - \sqrt{m} f\left(\frac{j}{T}\right) - \sqrt{m} \eta_j\right\}.$$

This is the deterministic component of the approximation error due to binning. It is easy to see that

$$\min_{(j-1)m+1 \leq i \leq jm} \left[f\left(\frac{i}{n}\right) - f\left(\frac{j}{T}\right)\right] \tag{29}$$
$$\leq X_j - \eta_j - f\left(\frac{j}{T}\right) \leq \max_{(j-1)m+1 \leq i \leq jm} \left[f\left(\frac{i}{n}\right) - f\left(\frac{j}{T}\right)\right].$$

Since $f$ is in a Hölder ball with smoothness $d = \min(\alpha - \frac{1}{p}, 1)$, then equation (29) implies

$$
\begin{aligned}
|\epsilon_j| &\leq \sqrt{m} E \left| X_j - f\left(\frac{j}{T}\right) - \eta_j \right| \\
&\leq \sqrt{m} \max_{(j-1)m+1 \leq i \leq jm} \left| f\left(\frac{i}{n}\right) - f\left(\frac{j}{T}\right) \right| \leq C\sqrt{m} T^{-d}.
\end{aligned}
\tag{30}
$$

Set

$$\zeta_j = \sqrt{m} X_j - \sqrt{m} f\left(\frac{j}{T}\right) - \sqrt{m} b_m - \epsilon_j - \frac{1}{2} Z_j.$$

Then $E\zeta_j = 0$ and $\sqrt{m} X_j = \sqrt{m} f(j/T) + \epsilon_j + \frac{1}{2} Z_j + \zeta_j$. The random error $\zeta_j$ is the sum of two terms, $\zeta_{1j} = \sqrt{m} X_j - \sqrt{m} f(j/T) - \sqrt{m} \eta_j - \epsilon_j$ and $\zeta_{2j} = \sqrt{m} \eta_j - \frac{1}{2} Z_j$, where $\zeta_{1j}$ is the random component of the approximation error due to binning, and $\zeta_{2j}$ is the error of approximating the median by the Gaussian variable. From equation (29) we have $|\zeta_{1j}| \leq C\sqrt{m} T^{-d}$ and so

$$E|\zeta_{1j}|^l \leq C_l m^{l/2} T^{-dl}. \tag{31}$$

A bound for the approximation error $\zeta_{2j}$ is given in Corollary 1,

$$|\zeta_{2j}| \leq \frac{C}{m^{1/2}}(1 + |Z_j|^2) \qquad \text{when } |Z_j| \leq \varepsilon\sqrt{m} \tag{32}$$

for some $\varepsilon > 0$, and the probability of $|Z_j| > \varepsilon\sqrt{m}$ is exponentially small. Hence for any finite integer $l \geq 1$ (here $l$ is fixed and $m = n^\gamma \to \infty$),

$$
\begin{aligned}
E|\zeta_{2j}|^l &= E|\zeta_{2j}|^l \{|Z_j| \leq \varepsilon\sqrt{m}\} + E|\zeta_{2j}|^l \{|Z_j| > \varepsilon\sqrt{m}\} \\
&\leq C_l m^{-l/2} + (E|\zeta_{2j}|^{2l})^{1/2} [P\{|Z_j| > \varepsilon\sqrt{m}\}]^{1/2}
\end{aligned}
$$

for some constant $C_l > 0$, where

$$P\{|Z| > \varepsilon\sqrt{m}\} \leq \frac{1}{2} \exp\left(-\frac{\varepsilon^2}{2} m\right)$$

by Mill's ratio inequality

$$\frac{\varphi(x)}{1 - \Phi(x)} > \max\left\{x, \frac{2}{\sqrt{2\pi}}\right\} \geq \frac{1}{2}\left(x + \frac{2}{\sqrt{2\pi}}\right) \qquad \text{for } x > 0 \tag{33}$$

and

$$E|\sqrt{m}\eta_j|^{2l} \leq m^l E|\eta_j|^{2l} \leq D_l m^l \tag{34}$$

for some constant $D_l > 0$ because of Assumption (A2), so we have

$$E|\zeta_{2j}|^l \leq C_l m^{-l/2}. \tag{35}$$

Details for equation (34) are as follows. Assumption (A2) implies

$$P(|\xi_i| \geq |x|) \leq \frac{C}{|x|^{\epsilon_3}}.$$

For $m = 2v + 1$ i.i.d. $\xi_i$, from equation (65) the density of the sample median is

$$g(x) = \frac{\sqrt{8v}}{\sqrt{2\pi}} [4H(x)(1 - H(x))]^v h(x) \exp\left(O\left(\frac{1}{v}\right)\right)$$

$$\leq \frac{\sqrt{8v}}{\sqrt{2\pi}} \left[\frac{4C}{|x|^{\epsilon_3}}\right]^v h(x) \exp\left(O\left(\frac{1}{v}\right)\right)$$

$$= \frac{\sqrt{8v}}{\sqrt{2\pi}} \left[\frac{4C}{|x|^{\epsilon_3/2}}\right]^v \frac{1}{|x|^{v\epsilon_3/2}} h(x) \exp\left(O\left(\frac{1}{v}\right)\right).$$

When $|x|^{\epsilon_3/2} \geq 8C$, we have

$$\frac{\sqrt{8v}}{\sqrt{2\pi}} \left[\frac{4C}{|x|^{\epsilon_3/2}}\right]^v \leq \frac{\sqrt{8v}}{\sqrt{2\pi} 2^v}$$

which is bounded for all $v$. This implies as $v \to \infty$ ($m \sim n^\gamma$ in our procedure) the median has any finite moments.

Thus we have

$$E|\zeta_j|^l \leq 2^{l-1}(E|\zeta_{1j}|^l + E|\zeta_{2j}|^l) \leq C_l m^{-l/2} + C_l m^{l/2} T^{-dl}$$

from equations (31) and (35). Equation (28) then follows from Chebyshev's inequality. □

REMARK 5. In the proof of Proposition 2, we will see that the noise $\zeta_j$ has negligible contribution to the risk of our procedure comparing with the Gaussian noise $\frac{1}{2}Z_j$, when the tail bound $P(|\zeta_j| > a)$ decays faster than any polynomial of $n$. For $m = n^\gamma$ we have $T^{2d}/m = n^{2d-\gamma(2d+1)}$. Then from equation (28) it is enough to require $0 < \gamma < \frac{2d}{2d+1}$, that is,

(36) $$d = \min\left(\alpha - \frac{1}{p}, 1\right) > \frac{\gamma}{2(1-\gamma)}$$

which is satisfied under our assumption (see also Remark 7).

REMARK 6. In the proofs of Theorems 3 and 4, we shall assume without loss of generality that $h(0)$ is known and equal to 1 since it can be estimated accurately in the sense that there is an estimator $\widehat{h}(0)$ such that

(37) $$P\{|\widehat{h}^{-2}(0) - h^{-2}(0)| > n^{-\delta}\} \leq c_l n^{-l}$$

for some $\delta > 0$ and all $l \geq 1$. For instance, we may estimate $h^{-2}(0)$ by

(38) $$\widehat{h}^{-2}(0) = \frac{8m}{T} \sum (X_{2k-1} - X_{2k})^2.$$

Note that $E \frac{m}{T/2} \sum (X_{2k-1} - X_{2k})^2 = \frac{1}{4} h^{-2}(0) + O(\sqrt{m} T^{-d})$, and it is easy to show

$$E \left| \frac{8m}{T} \sum (X_{2k-1} - X_{2k})^2 - h^{-2}(0) \right|^l \leq C_l (\sqrt{m} T^{-d})^l$$

where $\sqrt{m} T^{-d} = n^{-\delta}$ with $\delta > 0$ in our assumption. Then equation (37) holds by Chebyshev inequality. It is very important to see that the asymptotic risk properties of our estimator (16) does not change when replacing $\lambda_*$ by $\lambda_*(1 + O(n^{-\delta}))$, thus in the rest of our analysis we may just assume that $h(0) = 1$ without loss of generality.

We now consider the wavelet transform of the medians of the binned data. From Proposition 1 we may write

$$\frac{1}{\sqrt{T}} X_i = \frac{g(i/T)}{\sqrt{T}} + \frac{\epsilon_i}{\sqrt{n}} + \frac{Z_i}{2\sqrt{n}} + \frac{\zeta_i}{\sqrt{n}}.$$

Let $(y_{j,k}) = T^{-1/2} W \cdot X$ be the discrete wavelet transform of the binned data. Then one may write

(39) $$y_{j,k} = \theta'_{j,k} + \epsilon_{j,k} + \frac{1}{2\sqrt{n}} z_{j,k} + \xi_{j,k}$$

where $\theta'_{j,k}$ are the discrete wavelet transform of $(g(\frac{i}{T}))_{1 \leq i \leq T}$ which are approximately equal to the true wavelet coefficients of $g$, $z_{j,k}$ are the transform of the $Z_i$'s and so are i.i.d. $N(0,1)$ and $\epsilon_{j,k}$ and $\xi_{j,k}$ are respectively the transforms of $(\frac{\epsilon_i}{\sqrt{n}})$ and $(\frac{\zeta_i}{\sqrt{n}})$. The following proposition studies the risk of BlockJS procedure in Step 2. For each single block the risk bounds here for BlockJS are similar to results in Cai (1999) where Gaussian noise was considered. But in the current setting the error terms $\epsilon_{j,k}$ and $\xi_{j,k}$ make the problem more complicated.

PROPOSITION 2. *Let the empirical wavelet coefficients $y_{j,k} = \theta'_{j,k} + \epsilon_{j,k} + \frac{1}{2\sqrt{n}} z_{j,k} + \xi_{j,k}$ be given as in (39) and let the block thresholding estimator $\hat{\theta}_{j,k}$ be defined as in (16). Then:*

(i) *for some constant $C > 0$*

(40) $$E \sum_{(j,k) \in B_j^i} (\hat{\theta}_{j,k} - \theta'_{j,k})^2 \leq \min\left\{ 4 \sum_{(j,k) \in B_j^i} (\theta'_{j,k})^2, 8\lambda_* L n^{-1} \right\} + 6 \sum_{(j,k) \in B_j^i} \epsilon_{j,k}^2 + C L n^{-2};$$

(ii) *for any $0 < \tau < 1$, there exists a constant $C_\tau > 0$ depending on $\tau$ only such that for all $(j,k) \in B_j^i$*

$$(41) \quad E(\hat{\theta}_{j,k} - \theta'_{j,k})^2 \leq C_\tau \cdot \min\left\{\max_{(j,k) \in B_j^i}\{(\theta'_{j,k} + \epsilon_{j,k})^2\}, Ln^{-1}\right\} + n^{-2+\tau}.$$

We need the following lemmas to prove Proposition 2. These three lemmas are from Brown et al. (2006). See also Cai (1999).

LEMMA 1. *Let $X_1, \ldots, X_n$ be independent random variables with $E(X_i) = 0$ for $i = 1, \ldots, n$. Suppose that $E|X_i|^k < M_k$ for all $i$ and all $k > 0$ with $M_k > 0$ some constant not depending on $n$. Let $Y = WX$ be an orthogonal transform of $X = (X_1, \ldots, X_n)'$. Then there exist constants $M'_k$ not depending on $n$ such that $E|Y_i|^k < M'_k$ for all $i = 1, \ldots, n$ and all $k > 0$.*

LEMMA 2. *Suppose $y_i = \theta_i + z_i, i = 1, \ldots, L$, where $\theta_i$ are constants and $z_i$ are random variables. Let $S^2 = \sum_{i=1}^L y_i^2$ and let $\hat{\theta}_i = (1 - \frac{\lambda L}{S^2})_+ y_i$. Then*

$$(42) \quad E\|\hat{\theta} - \theta\|_2^2 \leq \|\theta\|_2^2 \wedge 4\lambda L + 4E[\|z\|_2^2 I(\|z\|_2^2 > \lambda L)].$$

LEMMA 3. *Let $X \sim \chi_L^2$ and $\lambda > 1$. Then*

$$(43) \quad \begin{aligned} P(X \geq \lambda L) &\leq e^{-(L/2)(\lambda - \log \lambda - 1)} \quad and \\ EXI(X \geq \lambda L) &\leq \lambda L e^{-(L/2)(\lambda - \log \lambda - 1)}. \end{aligned}$$

PROOF OF PROPOSITION 2. We only give the proof for (i). From Proposition 1, we have $|\epsilon_j| \leq C\sqrt{m}T^{-d}$ and $\epsilon_{j,k} = \sum_i \frac{\epsilon_i}{\sqrt{n}} \int \phi_{J,i} \psi_{j,k}$. Hence

$$(44) \quad |\epsilon_{j,k}| \leq \sup_x \sum_i \left|\frac{\epsilon_i}{\sqrt{n}} \phi_{J,i}(x)\right| \cdot \int |\psi_{j,k}(x)|\,dx \leq CT^{-d}2^{-j/2}.$$

This, as well as Proposition 1, yields that

$$(45) \quad \sum_j \sum_k \epsilon_{j,k}^2 = \frac{1}{n}\sum_i \epsilon_i^2 \leq CT^{-2d}.$$

It is easy to see from Lemma 1 and Proposition 1 that

$$(46) \quad E|\xi_{j,k}|^l \leq C'_l(mn)^{-l/2} + C'_l(T^{2d}n/m)^{-l/2}$$

and for any $a > 0$

$$(47) \quad P(|\xi_{j,k}| > a) \leq C'_l(a^2mn)^{-l/2} + C'_l(a^2 T^{2d}n/m)^{-l/2}.$$

It follows from Lemma 2 that

$$E \sum_{(j,k) \in B_j^i} (\hat{\theta}_{j,k} - \theta'_{j,k})^2$$

$$\leq 2E \sum_{(j,k) \in B_j^i} [\hat{\theta}_{j,k} - (\theta'_{j,k} + \epsilon_{j,k})]^2 + 2 \sum_{(j,k) \in B_j^i} \epsilon_{j,k}^2$$

$$\leq 2\min\left\{ \sum_{(j,k) \in B_j^i} (\theta'_{j,k} + \epsilon_{j,k})^2, 4\lambda_* L n^{-1} \right\} + 2 \sum_{(j,k) \in B_j^i} \epsilon_{j,k}^2$$

$$+ 8E \sum_{(j,k) \in B_j^i} \left(\frac{1}{2\sqrt{n}} z_{j,k} + \xi_{j,k}\right)^2 I\left( \sum_{(j,k) \in B_j^i} \left(\frac{1}{2\sqrt{n}} z_{j,k} + \xi_{j,k}\right)^2 > \frac{\lambda_* L}{4n} \right)$$

$$\leq \min\left\{ 4 \sum_{(j,k) \in B_j^i} (\theta'_{j,k})^2, 8\lambda_* L n^{-1} \right\} + 6 \sum_{(j,k) \in B_j^i} \epsilon_{j,k}^2$$

$$+ 2n^{-1} E \sum_{(j,k) \in B_j^i} (z_{j,k} + 2\sqrt{n}\xi_{j,k})^2 I\left( \sum_{(j,k) \in B_j^i} (z_{j,k} + 2\sqrt{n}\xi_{j,k})^2 > \lambda_* L \right).$$

Denote by $A$ the event that all $|\xi_{j,k}|$ are bounded by $\frac{1}{2\sqrt{n}L}$, that is

$$A = \{|2\sqrt{n}\xi_{j,k}| \leq L^{-1} \text{ for all } (j,k) \in B_j^i\}.$$

Then it follows from (47) that for any $l \geq 1$

$$P(A^c) \leq \sum_{(j,k) \in B_j^i} P(|2\sqrt{n}\xi_{j,k}| > L^{-1})$$

(48)
$$\leq C'_l (L^{-2} m)^{-l/2} + C'_l (L^{-2} T^d / m)^{-l/2}.$$

Hence

$$D = E \sum_{(j,k) \in B_j^i} (z_{j,k} + 2\sqrt{n}\xi_{j,k})^2 I\left( \sum_{(j,k) \in B_j^i} (z_{j,k} + 2\sqrt{n}\xi_{j,k})^2 > \lambda_* L \right)$$

$$= E \sum_{(j,k) \in B_j^i} (z_{j,k} + 2\sqrt{n}\xi_{j,k})^2 I\left( A \cap \sum_{(j,k) \in B_j^i} (z_{j,k} + 2\sqrt{n}\xi_{j,k})^2 > \lambda_* L \right)$$

$$+ E \sum_{(j,k) \in B_j^i} (z_{j,k} + 2\sqrt{n}\xi_{j,k})^2 I\left( A^c \cap \sum_{(j,k) \in B_j^i} (z_{j,k} + 2\sqrt{n}\xi_{j,k})^2 > \lambda_* L \right)$$

$$= D_1 + D_2.$$

Note that for any $L > 1$, $(x+y)^2 \leq \frac{L}{L-1}x^2 + Ly^2$ for all $x$ and $y$. It then follows from Lemma 3 and Hölder's inequality that

$$D_1 = E \sum_{(j,k) \in B_j^i} (z_{j,k} + 2\sqrt{n}\xi_{j,k})^2 I\left(A \cap \sum_{(j,k) \in B_j^i} (z_{j,k} + 2\sqrt{n}\xi_{j,k})^2 > \lambda_* L\right)$$

$$\leq 2E \sum_{(j,k) \in B_j^i} z_{j,k}^2 I\left(\sum_{(j,k) \in B_j^i} z_{j,k}^2 > \lambda_* L - \lambda_* - 1\right)$$

$$+ 8nE \sum_{(j,k) \in B_j^i} \xi_{j,k}^2 I\left(\sum_{(j,k) \in B_j^i} z_{j,k}^2 > \lambda_* L - \lambda_* - 1\right)$$

$$\leq 2(\lambda_* L - \lambda_* - 1)e^{-L/2(\lambda_* - (\lambda_*+1)L^{-1} - \log(\lambda_* - (\lambda_*+1)L^{-1}) - 1)}$$

$$+ 8n \sum_{(j,k) \in B_j^i} (E\xi_{j,k}^{2p})^{1/p} \left(P\left(\sum_{(j,k) \in B_j^i} z_{j,k}^2 > \lambda_* L - \lambda_* - 1\right)\right)^{1/q},$$

where $p, q > 1$ and $\frac{1}{p} + \frac{1}{q} = 1$. For $m = n^\epsilon$ we take $\frac{1}{q} = 1 - \epsilon$. Then it follows from Lemma 3 and (46) that

$$D_1 \leq \lambda_* e^{(\lambda_*+1)/2} Ln^{-1} + CLm^{-1}n^{-1-\epsilon} = CLn^{-1}.$$

On the other hand, it follows from (46) and (48) (by taking $l$ sufficiently large) that

$$D_2 = E \sum_{(j,k) \in B_j^i} (z_{j,k} + 2\sqrt{n}\xi_{j,k})^2 I\left(A^c \cap \sum_{(j,k) \in B_j^i} (z_{j,k} + 2\sqrt{n}\xi_{j,k})^2 > \lambda_* L\right)$$

$$\leq E \sum_{(j,k) \in B_j^i} (2z_{j,k}^2 + 8n\xi_{j,k}^2) I(A^c)$$

$$\leq \sum_{(j,k) \in B_j^i} [2(Ez_{j,k}^4)^{1/2} + 8n(E\xi_{j,k}^4)^{1/2}] \cdot (P(A^c))^{1/2}$$

$$\leq n^{-1}.$$

Hence, $D = D_1 + D_2 \leq CLn^{-1}$ and consequently

$$E \sum_{(j,k) \in B_j^i} (\hat{\theta}_{j,k} - \theta'_{j,k})^2 \leq \min\left\{4 \sum_{(j,k) \in B_j^i} (\theta'_{j,k})^2, 8\lambda_* Ln^{-1}\right\}$$

$$+ 6 \sum_{(j,k) \in B_j^i} \epsilon_{j,k}^2 + CLn^{-2}$$

for some constant $C > 0$. □

Recall that $\theta'_{j,k}$'s are the discrete wavelet transform of $(f(\frac{i}{T}))_{1 \leq i \leq T}$ and $\theta_{j,k}$'s are true wavelet coefficients of $f$. The following lemma will be used to bound the difference of $\theta'_{j,k}$'s and $\theta_{j,k}$'s. The proof is straightforward and is thus omitted.

LEMMA 4. *Let $T = 2^J$ and let $f_J(x) = \sum_{k=1}^{T} \frac{1}{\sqrt{T}} f(\frac{k}{T}) \phi_{J,k}(x)$. Then*

$$\sup_{f \in B^\alpha_{p,q}(M)} \|f_J - f\|_2^2 \leq CT^{-2d} \qquad \text{where } d = \min(\alpha - 1/p, 1).$$

*Also, $|\theta'_{j,k} - \theta_{j,k}| \leq CT^{-d} 2^{-j/2}$ and consequently $\sum_{j=j_0}^{J-1} \sum_k (\theta'_{j,k} - \theta_{j,k})^2 \leq CT^{-2d}$.*

LEMMA 5. *Let $b_m$ and $\hat{b}_m$ be defined as in (9) and (17), respectively. Then*

$$\sup_{h \in \mathcal{H}} \left| b_m + \frac{h'(0)}{8h^3(0)m} \right| \leq Cm^{-2}, \tag{49}$$

$$\sup_{h \in \mathcal{H}} \sup_{f \in B^\alpha_{p,q}(M)} E(\hat{b}_m - b_m)^2 \leq C \cdot \max\{T^{-2d}, m^{-4}\}. \tag{50}$$

PROOF. It suffices to consider the case that $m = 2v + 1$ with $v \in \mathbb{N}$ (cf. Remark 1), then

$$E\xi_{\text{med}} = \int x \frac{(2v+1)!}{(v!)^2} H^v(x)(1 - H(x))^v \, dH(x),$$

where $H$ is the distribution function of $\xi_1$. For any $\delta > 0$, set $A_\delta = \{x : |H(x) - \frac{1}{2}| \leq \delta\}$. It follows from the definition of $\mathcal{H}$ that there exists a constant $\delta > 0$ such that for some $\epsilon > 0$ we have

$$|h^{(3)}(x)| \leq 1/\epsilon \quad \text{and} \quad \epsilon \leq h(x) \leq 1/\epsilon \tag{51}$$

uniformly over all $h \in \mathcal{H}$ for all $x \in A_\delta$. This property implies $H^{-1}(x)$ is well defined and differentiable up to the fourth order for $x \in A_\delta$. Decompose the expectation of the median into two parts:

$$E\xi_{\text{med}} = \left( \int_{A_\delta} + \int_{A_\delta^c} \right) x \frac{(2v+1)!}{(v!)^2} H^v(x)(1 - H(x))^v \, dH(x) \equiv Q_1 + Q_2.$$

Since the median has finite moments from equation (34), it is easy to see $Q_2$ decays to 0 exponentially fast as $v = O(n^{1/4}) \to \infty$ by the Cauchy–Schwarz inequality and tail probability equations (63) and (64). We now turn to $Q_1$. Note that

$$Q_1 = \int_{1/2-\delta}^{1/2+\delta} \left(H^{-1}(x) - H^{-1}\left(\frac{1}{2}\right)\right) \frac{(2v+1)!}{(v!)^2} x^v (1-x)^v \, dx$$

$$= \int_{1/2-\delta}^{1/2+\delta} \left[\frac{1}{2}(H^{-1})''\left(\frac{1}{2}\right)\left(x-\frac{1}{2}\right)^2 + \frac{(H^{-1})^{(4)}(\varsigma)}{24}\left(x-\frac{1}{2}\right)^4\right]$$

$$\times \frac{(2v+1)!}{(v!)^2} x^v (1-x)^v \, dx$$

since $x^v(1-x)^v$ is symmetric around $x = \frac{1}{2}$. Note that $\frac{(2v+1)!}{(v!)^2} x^v (1-x)^v$ is the density function of Beta$(v+1, v+1)$, and equation (51) implies that $(H^{-1})^{(4)}(\varsigma)$ is uniformly bounded over all $h \in \mathcal{H}$, then

$$Q_1 = \frac{1}{2}(H^{-1})''\left(\frac{1}{2}\right)\frac{(v+1)^2}{(2v+2)^2(2v+3)} + O\left(\frac{1}{m^2}\right) = -\frac{h'(0)}{8h^3(0)m} + O\left(\frac{1}{m^2}\right)$$

and (49) is established.

Note that for $m = 2v + 1$, $\lfloor \frac{m}{2} \rfloor = v$. From Proposition 1 we have

$$X_j = f\left(\frac{j}{T}\right) + b_m + \frac{1}{2\sqrt{m}} Z_j + \frac{1}{\sqrt{m}} \epsilon_j + \frac{1}{\sqrt{m}} \zeta_j.$$

Similarly we may write

$$X_j^* = f\left(\frac{j-1/2}{T}\right) + b_v + \frac{1}{2\sqrt{v}} Z_j^* + \frac{1}{\sqrt{v}} \epsilon_j^* + \frac{1}{\sqrt{v}} \zeta_j^*$$

with $Z_j^*$, $\epsilon_j^*$ and $\zeta_j^*$ satisfying properties (i), (ii), (iii) of Proposition 1, respectively. Then $\hat{b}_m - b_m = \frac{1}{T} \sum_j (X_j^* - X_j) - b_m$ can be written as a sum of five terms as follows:

$$\hat{b}_m - b_m = \frac{1}{T} \sum_j \left(f\left(\frac{j-1/2}{T}\right) - f\left(\frac{j}{T}\right)\right) + (b_v - 2b_m)$$

$$+ \left[\frac{1}{\sqrt{v}}\frac{1}{T}\sum_j \epsilon_j^* - \frac{1}{\sqrt{m}}\frac{1}{T}\sum_j \epsilon_j\right]$$

$$+ \left[\frac{1}{2\sqrt{v}}\frac{1}{T}\sum_j Z_j^* - \frac{1}{2\sqrt{m}}\frac{1}{T}\sum_j Z_j\right]$$

$$+ \left[\frac{1}{\sqrt{v}}\frac{1}{T}\sum_j \zeta_j^* - \frac{1}{\sqrt{m}}\frac{1}{T}\sum_j \zeta_j\right]$$

$$\equiv R_1 + R_2 + R_3 + R_4 + R_5.$$

It is easy to see that $\sup_{f\in B^\alpha_{p,q}(M)} R_1^2 \leq CT^{-2d}$ and $\sup_{h\in\mathcal{H}} R_2^2 \leq Cm^{-4}$. Proposition 1 yields $\sup_{h\in\mathcal{H}, f\in B^\alpha_{p,q}(M)} R_3^2 \leq CT^{-2d}$. Note that $Z_j^* - Z_j$ are independent for $j=1,\ldots,T$. So $ER_4^2 \leq \frac{1}{h^2(0)}(\frac{1}{v}+\frac{1}{m})\frac{1}{T} \leq Cn^{-1}$. Similarly, $\zeta_j^* - \zeta_j$ are independent and it then follows from Proposition 1 that $ER_5^2 = o(n^{-1})$. Hence,

$$\sup_{h\in\mathcal{H}, f\in B^\alpha_{p,q}(M)} E(\hat{b}_m - b_m)^2 \leq 5R_1^2 + 5R_2^2 + 5R_3^2 + 5ER_4^2 + 5ER_5^2$$

$$\leq C\max\{T^{-2d}, m^{-4}\}. \qquad \square$$

6.2. *Global adaptation: Proof of Theorem 3.* Let $\hat{f}_n$ be given as in (18). Note that

$$E\|\hat{f}_n - f\|_2^2 \leq 2E\|\hat{g}_n - g\|_2^2 + 2E(\hat{b}_m - b_m)^2.$$

Lemma 5 yields that $E(\hat{b}_m - b_m)^2 = o(n^{-2\alpha/(2\alpha+1)})$ and so we need only to focus on bounding $E\|\hat{g}_n - g\|_2^2$. Note that the functions $f$ and $g$ differ only by a constant $b_m$ and so the wavelet coefficients coincide, that is, $\theta_{j,k} = \int f\psi_{j,k} = \int g\psi_{j,k}$. Decompose $E\|\hat{g}_n - g\|_2^2$ into three terms as follows:

$$E\|\hat{g}_n - g\|_2^2 = \sum_k E(\hat{\tilde{\theta}}_{j_0,k} - \tilde{\theta}_{j,k})^2 + \sum_{j=j_0}^{J-1}\sum_k E(\hat{\theta}_{j,k} - \theta_{j,k})^2 + \sum_{j=J}^{\infty}\sum_k \theta_{j,k}^2$$

(52)
$$\equiv S_1 + S_2 + S_3.$$

It is easy to see that the first term $S_1$ and the third term $S_3$ are small:

(53) $$S_1 = 2^{j_0} n^{-1} \epsilon^2 = o(n^{-2\alpha/(1+2\alpha)}).$$

Note that for $x\in\mathbb{R}^m$ and $0 < p_1 \leq p_2 \leq \infty$

(54) $$\|x\|_{p_2} \leq \|x\|_{p_1} \leq m^{1/p_1 - 1/p_2}\|x\|_{p_2}.$$

Since $f\in B^\alpha_{p,q}(M)$, so $2^{js}(\sum_{k=1}^{2^j}|\theta_{j,k}|^p)^{1/p} \leq M$. Now (54) yields that

(55) $$S_3 = \sum_{j=J}^{\infty}\sum_k \theta_{j,k}^2 \leq C2^{-2J(\alpha\wedge(\alpha+1/2-1/p))}.$$

Proposition 2, Lemma 4 and equation (45) yield that

$$S_2 \leq 2\sum_{j=j_0}^{J-1}\sum_k E(\hat{\theta}_{j,k} - \theta'_{j,k})^2 + 2\sum_{j=j_0}^{J-1}\sum_k (\theta'_{j,k} - \theta_{j,k})^2$$

(56)
$$\leq \sum_{j=j_0}^{J-1}\sum_{i=1}^{2^j/L}\min\left\{8\sum_{(j,k)\in B_j^i}\theta_{j,k}^2, 8\lambda_* Ln^{-1}\right\} + 6\sum_{j=j_0}^{J-1}\sum_k \epsilon_{j,k}^2$$

ROBUST NONPARAMETRIC ESTIMATION 25

$$+ Cn^{-1} + 10 \sum_{j=j_0}^{J-1} \sum_k (\theta'_{j,k} - \theta_{j,k})^2$$

$$\leq \sum_{j=j_0}^{J-1} \sum_{i=1}^{2^j/L} \min\left\{ 8 \sum_{(j,k) \in B_j^i} \theta_{j,k}^2, 8\lambda_* L n^{-1} \right\} + Cn^{-1} + CT^{-2d}.$$

We now divide into two cases. First consider the case $p \geq 2$. Let $J_1 = [\frac{1}{1+2\alpha} \times \log_2 n]$. So, $2^{J_1} \approx n^{1/(1+2\alpha)}$. Then (56) and (54) yield

$$S_2 \leq 8\lambda_* \sum_{j=j_0}^{J_1-1} \sum_{i=1}^{2^j/L} Ln^{-1} + 8 \sum_{j=J_1}^{J-1} \sum_k \theta_{j,k}^2 + Cn^{-1} + CT^{-2d} \leq Cn^{-2\alpha/(1+2\alpha)}.$$

By combining this with (53) and (55), we have $E\|\hat{f}_n - f\|_2^2 \leq Cn^{-2\alpha/(1+2\alpha)}$ for $p \geq 2$.

Now let us consider the case $p < 2$. First we state the following lemma without proof.

LEMMA 6. *Let $0 < p < 1$ and $S = \{x \in I\mathbb{R}^k : \sum_{i=1}^k x_i^p \leq B, x_i \geq 0, i = 1, \ldots, k\}$. Then for $A > 0$, $\sup_{x \in S} \sum_{i=1}^k (x_i \wedge A) \leq B \cdot A^{1-p}$.*

Let $J_2$ be an integer satisfying $2^{J_2} \asymp n^{1/(1+2\alpha)}(\log n)^{(2-p)/p(1+2\alpha)}$. Note that

$$\sum_{i=1}^{2^j/L} \left( \sum_{(j,k) \in B_j^i} \theta_{j,k}^2 \right)^{p/2} \leq \sum_{k=1}^{2^j} (\theta_{j,k}^2)^{p/2} \leq M 2^{-jsp}.$$

It then follows from Lemma 6 that

(57)
$$\sum_{j=J_2}^{J-1} \sum_{i=1}^{2^j/L} \min\left\{ 8 \sum_{(j,k) \in B_j^i} \theta_{j,k}^2, 8\lambda_* L n^{-1} \right\}$$
$$\leq Cn^{-2\alpha/(1+2\alpha)}(\log n)^{(2-p)/(p(1+2\alpha))}.$$

On the other hand,

(58)
$$\sum_{j=j_0}^{J_2-1} \sum_{i=1}^{2^j/L} \min\left\{ 8 \sum_{(j,k) \in B_j^i} \theta_{j,k}^2, 8\lambda_* L n^{-1} \right\}$$
$$\leq \sum_{j=j_0}^{J_2-1} \sum_b 8\lambda_* L n^{-1} \leq Cn^{-2\alpha/(1+2\alpha)}(\log n)^{(2-p)/(p(1+2\alpha))}.$$

We finish the proof for the case $p < 2$ by putting (53), (55), (57) and (58) together:

$$E\|\hat{f}_n - f\|_2^2 \leq Cn^{-2\alpha(1+2\alpha)}(\log n)^{(2-p)/(p(1+2\alpha))}.$$

REMARK 7. To make the risk of $\hat{b}_m$ negligible we need to have $m^{-4} = o(n^{-2\alpha/(1+2\alpha)})$ (see Lemma 5), and to make the approximation error $\|f_J - f\|_2^2$ negligible, we need to have $T^{-2((\alpha-1/p)\wedge 1)} = O(n^{-2\alpha/(1+2\alpha)})$ (see Lemma 4). These constraints lead to our choice of $m = n^{1/4}$ and $T = n^{3/4}$. Then we need $\frac{3}{2}(\alpha - \frac{1}{p}) > \frac{2\alpha}{1+2\alpha}$ or equivalently $\frac{2\alpha^2 - \alpha/3}{1+2\alpha} > \frac{1}{p}$. This last condition is purely due to approximation error over Besov spaces.

6.3. *Local adaptation: Proof of Theorem 4.* For simplicity, we give the proof for Hölder classes $\Lambda^\alpha(M)$ instead of local Hölder classes $\Lambda^\alpha(M, t_0, \delta)$. Note that for all $f \in \Lambda^\alpha(M)$, $|\theta_{j,k}| = |\langle f, \psi_{j,k}\rangle| \leq C2^{-j(1/2+\alpha)}$ for some constant $C > 0$ not depending on $f$. Note also that for any random variables $X_i$, $i = 1, \ldots, n$, $E(\sum_{i=1}^n X_i)^2 \leq (\sum_{i=1}^n (EX_i^2)^{1/2})^2$. It then follows that

$$E(\hat{f}_n(t_0) - f(t_0))^2$$

$$= E\left[\sum_{k=1}^{2^{j_0}}(\hat{\tilde{\theta}}_{j_0,k} - \tilde{\theta}_{j_0,k})\phi_{j_0,k}(t_0) + \sum_{j=j_0}^{\infty}\sum_{k=1}^{2^j}(\hat{\theta}_{j,k} - \theta_{j,k})\psi_{j,k}(t_0)\right.$$

$$\left. - (\hat{b}_m - b_m)\right]^2$$

$$\leq \left[(E(\hat{b}_m - b_m)^2)^{1/2} + \sum_{k=1}^{2^{j_0}}(E(\hat{\tilde{\theta}}_{j_0,k} - \tilde{\theta}_{j_0,k})^2\phi_{j_0,k}^2(t_0))^{1/2}\right.$$

$$\left. + \sum_{j=j_0}^{J-1}\sum_{k=1}^{2^j}(E(\hat{\theta}_{j,k} - \theta_{j,k})^2\psi_{j,k}^2(t_0))^{1/2} + \sum_{j=J}^{\infty}\sum_{k=1}^{2^j}|\theta_{j,k}\psi_{j,k}(t_0)|\right]^2$$

$$\equiv (Q_1 + Q_2 + Q_3 + Q_4)^2.$$

Lemma 5 yields that

(59) $$Q_1 = (E(\hat{b}_m - b_m)^2)^{1/2} = o(n^{-\alpha/(2\alpha+1)}).$$

Since the wavelet $\psi$ is compactly supported, so there are at most $N$ basis functions $\psi_{j,k}$ at each resolution level $j$ that are nonvanishing at $t_0$, where $N$ is the length of the support of $\psi$. Denote $K(t_0, j) = \{k : \psi_{j,k}(t_0) \neq 0\}$. Then $|K(t_0, j)| \leq N$. It is easy to see that both $Q_2$ and $Q_4$ are small:

(60) $$Q_2 = \sum_{k=1}^{2^{j_0}}(E(\hat{\tilde{\theta}}_{j_0,k} - \tilde{\theta}_{j_0,k})^2)^{1/2}|\phi_{j_0,k}(t_0)| = O(n^{-1})$$

and

$$(61) \quad Q_4 = \sum_{j=J}^{\infty} \sum_{k=1}^{2^j} |\theta_{j,k}||\psi_{j,k}(t_0)| \leq \sum_{j=J}^{\infty} N\|\psi\|_\infty 2^{j/2} C 2^{-j(1/2+\alpha)} \leq CT^{-\alpha}.$$

We now consider the third term $Q_3$. Applying the bound (41) in Proposition 2 with $\tau < 1/(1+2\alpha)$ together with Lemma 4 and the bound for $\epsilon_{j,k}$ given in (44), we have

$$Q_3 \leq \sum_{j=j_0}^{J-1} \sum_{k \in K(t_0,j)} 2^{j/2} \|\psi\|_\infty (E(\hat{\theta}_{j,k} - \theta_{j,k})^2)^{1/2}$$

$$(62) \quad \leq C \sum_{j=j_0}^{J-1} 2^{j/2} [\min(2^{-j(1+2\alpha)} + T^{-2(\alpha \wedge 1)} 2^{-j}, Ln^{-1}) + n^{-2+\tau}]^{\frac{1}{2}}$$

$$\leq C \left(\frac{\log n}{n}\right)^{\alpha/(1+2\alpha)}.$$

Combining equations (59)–(63) we have

$$E(\hat{f}_n(t_0) - f(t_0))^2 \leq C(\log n/n)^{2\alpha/(1+2\alpha)}.$$

6.4. *Proofs of Theorems 1 and 2.* Let $G(x)$ be the cumulative distribution function of $X_{\text{med}}$ and let $\varphi(z)$ and $\Phi(z)$ denote respectively the density and cumulative distribution function of a standard normal random variable. Using similar arguments in the proof of Theorem 3 in Zhou (2005) or a sketch in Section 6 of Komlós, Major and Tusnády (1975), we need only to show

$$(63) \quad G(x) = \Phi(\sqrt{8k}x) \exp(O(k|x|^3 + |x| + k^{-1/2})) \qquad \text{for } -\varepsilon \leq x \leq 0$$

and

$$(64) \quad 1 - G(x) = (1 - \Phi(\sqrt{8k}x)) \exp(O(k|x|^3 + |x| + k^{-1/2}))$$

$$\text{for } 0 \leq x \leq \varepsilon,$$

where $O(x)$ means a value between $-Cx$ and $Cx$ uniformly for some constant $C > 0$. Related asymptotic expansions for the distribution of median can be found in current literature, for instance, Burnashev (1996), but the major theorems there are not sufficient to establish the median coupling inequality.

Let $H(x)$ be distribution function of $X_1$. The density of the median $X_{(k+1)}$ is

$$g(x) = \frac{(2k+1)!}{(k!)^2} H^k(x)(1 - H(x))^k h(x).$$

Stirling's formula, $j! = \sqrt{2\pi} j^{j+1/2} \exp(-j + \epsilon_j)$ with $\epsilon_j = O(1/j)$, gives

$$g(x) = \frac{(2k+1)!}{4^k (k!)^2} [4H(x)(1-H(x))]^k h(x)$$

$$= \frac{2\sqrt{2k+1}}{e\sqrt{2\pi}} \left(\frac{2k+1}{2k}\right)^{2k+1} [4H(x)(1-H(x))]^k h(x) \exp\left(O\left(\frac{1}{k}\right)\right).$$

It is easy to see $|\sqrt{2k+1}/\sqrt{2k} - 1| \leq k^{-1}$, and

$$\left(\frac{2k+1}{2k}\right)^{2k+1} = \exp\left(-(2k+1)\log\left(1 - \frac{1}{2k+1}\right)\right) = \exp\left(1 + O\left(\frac{1}{k}\right)\right).$$

Then we have, when $0 < H(x) < 1$,

$$(65) \qquad g(x) = \frac{\sqrt{8k}}{\sqrt{2\pi}} [4H(x)(1-H(x))]^k h(x) \exp\left(O\left(\frac{1}{k}\right)\right).$$

From the Lipschitz assumption in the theorem, Taylor's expansion gives

$$4H(x)(1-H(x)) = 1 - 4(H(x) - H(0))^2$$

$$= 1 - 4\left[\int_0^x (h(t) - h(0))\,dt + h(0)x\right]^2$$

$$= 1 - 4(h(0)x + O(x^2))^2$$

for $0 \leq |x| \leq \varepsilon$, that is, $\log(4H(x)(1-H(x))) = -4h^2(0)x^2 + O(|x|^3)$ when $|x| \leq 2\varepsilon$ for some $\varepsilon > 0$. Here $\varepsilon$ is chosen small enough such that $h(x) > 0$ for $|x| \leq 2\varepsilon$. The Lipschitz assumption in the theorem also implies $\frac{h(x)}{h(0)} = 1 + O(|x|) = \exp(O(|x|))$ for $|x| \leq 2\varepsilon$. Thus

$$g(x) = \frac{\sqrt{8k}h(0)}{\sqrt{2\pi}} \exp(-8kh^2(0)x^2/2 + O(k|x|^3 + |x| + k^{-1})) \qquad \text{for } |x| \leq 2\varepsilon.$$

Now we approximate the distribution function of $X_{\text{med}}$ by the distribution function of normal random variable. Without loss of generality we assume $h(0) = 1$. We write

$$g(x) = \frac{\sqrt{8k}}{\sqrt{2\pi}} \exp(-8kx^2/2 + O(k|x|^3 + |x| + k^{-1})) \qquad \text{for } |x| \leq 2\varepsilon.$$

Now we use this approximation of density functions to give the desired approximation of distribution functions. Specifically we shall show

$$(66) \qquad G(x) = \int_{-\infty}^x g(t)\,dt \leq \Phi(\sqrt{8k}x) \exp(Ck|x|^3 + C|x| + Ck^{-1})$$

and

$$(67) \qquad G(x) \geq \Phi(\sqrt{8k}x) \exp(-Ck|x|^3 - C|x| - Ck^{-1})$$

for all $-\varepsilon \leq x \leq 0$ and some $C > 0$. The proof for $0 \leq x \leq \varepsilon$ is similar. Now we give the proof for inequality (66). Note that

$$
\begin{aligned}
(68) \quad &(\Phi(\sqrt{8k}x)\exp(-Ckx^3 - Cx + Ck^{-1}))' \\
&= \sqrt{8k}\varphi(\sqrt{8k}x)\exp(-Ckx^3 - Cx + Ck^{-1}) \\
&\quad - \Phi(\sqrt{8k}x)(3Ckx^2 - C)\exp(-Ckx^3 - Cx + Ck^{-1}).
\end{aligned}
$$

From Mill's ratio, inequality (33), we have $\Phi(\sqrt{8k}x)(-\sqrt{8k}x) < \varphi(\sqrt{8k}x)$ and hence

$$
\begin{aligned}
&-\Phi(\sqrt{8k}x)(3Ckx^2)\exp(-Ckx^3 - Cx + Ck^{-1}) \\
&\geq \sqrt{8k}\varphi(\sqrt{8k}x)\left(\frac{3C}{8}x\right)\exp(-Ckx^3 - Cx + Ck^{-1}).
\end{aligned}
$$

This and (68) yield

$$
\begin{aligned}
&(\Phi(\sqrt{8k}x)\exp(-Ckx^3 - Cx + Ck^{-1}))' \\
&\geq \sqrt{8k}\varphi(\sqrt{8k}x)\left(1 + \frac{3C}{8}x\right)\exp(-Ckx^3 - Cx + Ck^{-1}) \\
&\geq \sqrt{8k}\varphi(\sqrt{8k}x)\exp(Cx/2)\exp(-Ckx^3 - Cx + Ck^{-1}) \\
&\geq \sqrt{8k}\varphi(\sqrt{8k}x)\exp\left(-\frac{C}{2}kx^3 - \frac{C}{2}x + Ck^{-1}\right).
\end{aligned}
$$

Here in the second inequality we apply $(1 + C3x/8) \geq \exp(Cx/2)$ when $|Cx| \leq C(2\varepsilon) < 1/2$. Thus we have

$$
\begin{aligned}
&(\Phi(\sqrt{8k}x)\exp(-Ckx^3 - Cx + Ck^{-1}))' \\
&\geq \sqrt{8k}\varphi(\sqrt{8k}x)\exp(O(k|x|^3 + |x| + k^{-1}))
\end{aligned}
$$

for $C$ sufficiently large and for $-2\varepsilon \leq x \leq 0$, then

$$
\begin{aligned}
\int_{-2\varepsilon}^{x} g(t)\,dt &\leq \int_{-2\varepsilon}^{x} (\Phi(\sqrt{8k}t)\exp(-Ckt^3 - Ct + Ck^{-1}))' \\
&= \begin{bmatrix} \Phi(\sqrt{8k}x)\exp(-Ckx^3 - Cx + Ck^{-1}) \\ -\Phi(\sqrt{8k}\cdot(2\varepsilon))\exp(C(k(2\varepsilon)^3 + k^{-1})) \end{bmatrix} \\
&\leq \Phi(\sqrt{8k}x)\exp(-Ckx^3 - Cx + Ck^{-1}).
\end{aligned}
$$

In (65) we see

$$
\begin{aligned}
\int_{-\infty}^{-2\varepsilon} g(t)\,dt &= \int_{-\infty}^{-2\varepsilon} \frac{(2k+1)!}{(k!)^2} H^k(t)(1 - H(t))^k h(t)\,dt \\
&= \int_{0}^{H(-2\varepsilon)} \frac{(2k+1)!}{(k!)^2} u^k(1-u)^k\,du
\end{aligned}
$$

$$= o(k^{-1}) \int_{H(-3\varepsilon/2)}^{H(-\varepsilon)} \frac{(2k+1)!}{(k!)^2} u^k (1-u)^k \, du$$

$$\leq o(k^{-1}) \int_{H(-2\varepsilon)}^{H(x)} \frac{(2k+1)!}{(k!)^2} u^k (1-u)^k \, du$$

$$= o(k^{-1}) \int_{-2\varepsilon}^{x} g(t) \, dt,$$

where the third equality is from the fact that $u_1^k(1-u_1)^k = o(k^{-1})u_2^k(1-u_2)^k$ uniformly for $u_1 \in [0, H(-2\varepsilon)]$ and $u_2 \in [H(-3\varepsilon/2), H(-\varepsilon)]$. Thus we have

$$G(x) \leq \Phi(\sqrt{8k}x) \exp(-Ckx^3 - Cx + Ck^{-1}),$$

which is equation (66). Equation (67) can be established in a similar way.

REMARK. Note that in the proof of Theorem 1 it can be seen easily that constants $C$ and $\epsilon$ in equation (5) depend only on the ranges of $h(0)$ and the bound of Lipschitz constants of $h$ at a fixed open neighborhood of 0. Theorem 2 then follows from the proof of Theorem 1 together with this observation.

L. D. Brown
T. T. Cai
The Wharton School
University of Pennsylvania
Philadelphia, Pennsylvania 19104
USA
E-mail: lbrown@wharton.upenn.edu
          tcai@wharton.upenn.edu

H. H. Zhou
Department of Statistics
Yale University
New Haven, Connecticut 06511
USA
E-mail: huibin.zhou@yale.edu